\documentclass{article}
\usepackage{amssymb}
\usepackage{amsmath,amsthm,amsfonts}
\usepackage{graphicx}
\usepackage{enumerate}
\usepackage{latexsym}
\usepackage{comment}
\usepackage{engord}
\usepackage{color}
\usepackage{verbatim}
\newcommand{\ddx}[1]{\ensuremath{\frac{\partial}{\partial #1}}}
\newcommand{\inv}{\ensuremath{^{-1}}}
\newcommand{\real}{\ensuremath{\mathbb{R}}}
\newcommand{\ds}{\ensuremath{\displaystyle}}
\newcommand{\cplx}{\ensuremath{\mathbb{C}}}
\newcommand{\gothic}{\mathfrak}
\renewcommand{\gg}{{\gothic g}}
\newcommand{\ym}{{\mathfrak Y}{\gothic m}}

\newcommand{\fc}{{\mathcal F}_\cplx}

\newtheorem{thm}{Theorem}[section]
\newtheorem{cor}[thm]{Corollary}
\newtheorem{prop}[thm]{Proposition}
\newtheorem{dfn}[thm]{Definition}
\numberwithin{equation}{section}
\numberwithin{thm}{section}

\title{Geometric Aspects of the Kapustin-Witten Equations}
\author{ Michael Gagliardo$^1$ and Karen Uhlenbeck$^2$}

\begin{document}
\maketitle
\noindent    (1) Department of Mathematics, California Lutheran University, Thousand Oaks, California, USA, mgagliar@calllutheran.edu \\
    (2) Department of Mathematics, The University of Texas at Austin, Austin, Texas, USA, uhlen@math.utexs.edu\\

\tableofcontents

\begin{abstract}This expository article introduces the Kapustin-Witten equations to mathematicians.  We discuss the connections between the Complex Yang-Mills equations and the Kapustin-Witten equations.  In addition, we show the relation between the Kapustin-Witten equations, the moment map condition and the gradient Chern-Simons flow.  The new results in the paper correspond to estimates on the solutions to the KW equations given an estimate on the complex part of the connection. This leaves open the problem of obtaining global estimates on the complex part of the connection.\end{abstract}


{\bf Classification:} Differential Geometry (53), Global Analysis (58)

\pagebreak
The authors dedicate this paper to Richard Palais on the occasion of his 80\textsuperscript{th} birthday.  The senior author enthusiastically expresses her appreciation for the guidance she has received from Dick throughout her career, from graduate student days to the present.  In fact, Dick encouraged her to look at problems in gauge theory, which became a successful research experience. He serves as a wonderful model for a younger mathematician:  a successful research career,  a committed teacher to both undergraduate and graduate students, and a spectacularly second career in mathematical graphics software.

The junior author is greatly appreciative of the mentoring and motivation he has received during his, unfortunately few, interactions with Dick. The junior author feels very fortunate to be part of this academic family tree.

A brief outline of the paper follows.  We should start with the caveat that, with the exception of the material in the second part of section 4, this is an expository paper, following the spirit of some of Dick's best known work. We would like to introduce this subject to mathematicians.  Most of the results can be found embedded in the original paper of Kapustin and Witten \cite{kapustin2007}.

Chapter 1 is a review of the self-dual Yang-Mills equations \cite{MR518436}\cite{MR1081321}, and chapter 2 gives us the definition.  Chapter 3 goes into the relationship with the full complex Yang-Mills equations.  We delve in Chapter 4 into a little bit of analysis and give some estimates on solutions\cite{MR887284}\cite{springerlink:10.1007/BF01947069}\cite{MR648355}\cite{MR1809291}\cite{MR2030823}\cite{taubes2012}. Finally the geometry of the moment map and the flow for complex Chern-Simons are dealt with in Chapters 5 and 6\cite{witten2011}\cite{MR2039989}\cite{MR2359489}.

Happy Birthday Dick!

\section{Review of the Self-Dual Yang-Mills Equations}

The goal of this paper is to construct the complex Yang-Mills and Kapustin-Witten equations and to introduce basic properties. To begin, we will review the real Yang-Mills equations.  A more detailed introduction to the topic can be found in \cite{MR518436}.

The Yang-Mills equations are defined for a connection on a principle $G$-bundle, $P$, over a smooth 4-manifold $M$. Here $G$ is a compact gauge group, which may be set to $SU(2)$.  When $P$ is restricted to an open neighborhood ${\mathcal O}$ of $M$ it is locally isomorphic to ${\mathcal O}\times G$. We define a connection, $D_A=d+A$, on $P|_{\mathcal O}$, where $A$ (the gauge potential) is a $\gg$ valued one form on $M$ and $d$ is the standard differential on forms. In coordinates, we have
$$D_A=\sum_{j=1}^4\ddx{x^j}+A_jdx^j.$$
Here $A_j$ is a $\gg$ valued function on $M$.  From here we can define the curvature of the connection, $F_A$ (the gauge field) by
\begin{equation}F_A=D_A^2=[D_A,D_A]=\frac{1}{2}\sum_{j,k=1}^4\ddx{x^j}A_k-\ddx{x^k}A_j+\frac{1}{2}[A_j,A_k]dx^j\wedge dx^k.\end{equation}
On any Riemannian manifold we can define Yang-Mills functional on the space of connections by
\begin{equation}\ym(D_A)=\int_M|F_A|^2d\mu=\int_Mtr(F_A\wedge *F_A^*)\end{equation}
where $d\mu$ is the volume form on $M$ and $F_A^*=-F_A$.

The Yang-Mills equations are the Euler-Lagrange equations for the Yang-Mills functional and are defined in all dimensions. The equations are
\begin{equation}
(D^*_AF_A)_k=(*D_A*F_A)_k=[\nabla_j+A_j,F_{l,k}]\mu^{jl}=0
\end{equation}
where *, the covariant derivative $\nabla_j$ and the metric tensor $\mu^{ji}$ come from a Riemannian structure.  In this case, we note that 4 dimensions is special, since the equations $\ds D_A*F_A=0$ depend only on the conformal structure.

Recall from the definition of chern class, $c_2$, and Chern-Simons theory that
\begin{equation}\int_Mtr(F_A\wedge F_A)=\text{ topological term }+\int_{\partial M}CS(A).
\end{equation}
Here the topological term is $8\pi^2c_2$ when $\partial M=\varnothing$, and  $CS(A)$ is the Chern-Simons three-form when $\partial M\ne0$. In a flat bundle
\begin{equation}
CS(A)=A\wedge F_A-\frac{1}{6}(A\wedge A\wedge A).
\end{equation}

Finally, we have:

\begin{equation}\begin{aligned}
\ym(D_A)&=\int_{M}|F_A|^2d\mu=\frac{1}{2}\int_{M}|F_A\pm *F_A|^2d\mu\pm\int_{M}tr\left(F_A\wedge F_A\right)\\
&=2\int_{M}|F^\pm_A|^2d\mu\mp\int_{M}tr(F_A\wedge F_A)\\
&=2\int_{M}|F^\pm_A|^2d\mu\pm\text{ topological term }+\int_{\partial M}CS(A)
\end{aligned}\end{equation}

 We see that either equation $F^+_A=0$ or $F^-_A=0$ provides an absolute minimum in a fixed topological class with respect to fixed boundary data.  The study of solutions of the equations for $\partial M=\varnothing$ is the subject of Donaldson theory. For $\partial M\ne\varnothing$, in particular $M=M^3\times \real$, this is the subject of Floer theory.  In this case, $\ds F_A|_{M^3\times\pm\infty}=0$ and the self-dual connections $F^-_A=0$ provide flow lines between the flat connections.

\section{The Kapustin-Witten equations}

The Kapustin-Witten equations are defined on a Riemannian 4-manifold given a principle bundle $P$ with a (real) structure group $G$. For most present considerations, $G$ can be taken to be $SU(2)$ or $SO(3)$. The equations link a connection $D_A$ in $P$ with a one form $\phi$ which takes values in $P\times_G\gg$. The combination is naturally regarded as a complex connection $D_A+i\phi$. Here the configuration space $Q_\cplx$ is an affine complex vector space with complex conjugation defined.  The curvature, $\fc$, of the complex connection $D_A+i\phi$ is a two form with values in $P\times_G(\gg\otimes\cplx)$
\begin{equation}
\label{fc}\fc=[D_A+i\phi\wedge D_A+i\phi]=F_A-\frac{1}{2}[\phi\wedge\phi]+iD_A\phi.
\end{equation}
Here, $F_A$ is the curvature of the real connection $D_A$ and $D_A\phi$ is an extension of exterior differentiation, in coordinates
\begin{equation}
\label{zeroc}(D_A\phi)_{j,k}=\left(\left[\ddx{x^j}+A_j,\phi_k\right]-\left[\ddx{x^k}+A_k,\phi_j\right]\right)dx^j\wedge dx^k.
\end{equation}
Hence, complex flat connections satisfy an equation of the form
\begin{equation}
F_A-\frac{1}{2}[\phi\wedge\phi]=D_A\phi=0.
\end{equation}
These equations are not only invariant under the real gauge group $\mathcal{G}=C^\infty\left(P\times_G G\right)$, but are also invariant under the complex gauge group $\mathcal{G}_\cplx=C^\infty\left(P\times_G G_\cplx\right)$. Since the theory we seek is a theory with respect to the real gauge group, some method of choosing a point in the complex gauge orbit is necessary.

This is done by imposing the additional equation $D_A*\phi=0$. In section \ref{sec5}, we show this is the zero of the moment map for the real gauge group action using the trivial complex geometry we have introduced.  In two dimensions, the complex flat equations and the moment map equation
\begin{equation}\begin{aligned} F_A-\frac{1}{2}[\phi\wedge\phi]&=0,\\
D_A\phi&=0,\\
D_A*\phi&=0\
\end{aligned}
\end{equation}
form an equation which is elliptic modulo the real gauge group. These are Hitchin's equations and are the subject of much active research. However, in higher dimensions, these equations, while of great interest, are overdetermined.

In four dimensions, the system of of equations

$$\left(F_A-\frac{1}{2}[\phi\wedge\phi]\right)^+=0$$
$$\left(D_A\phi\right)^-=0$$
$$D_A*\phi=0$$
or alternatively
$$\left(F_A-\frac{1}{2}[\phi\wedge\phi]\right)^-=0$$
$$\left(D_A\phi\right)^+=0$$
$$D_A*\phi=0$$
form elliptic systems. However, these two systems fit into a family of equations

\begin{equation}\label{famKW}\begin{aligned}
\left(\cos{\theta}(F_A-\frac{1}{2}[\phi\wedge\phi])-\sin{\theta}D_A\phi\right)^+&=0\\
\left(\sin{\theta}(F_A-\frac{1}{2}[\phi\wedge\phi])+\cos{\theta}D_A\phi\right)^-&=0\\
D_A*\phi&=0
\end{aligned}\end{equation}
which are elliptic of the same index. This is the family of Kapustin-Witten equations. These can also be succinctly written as $e^{i\theta}\fc=*\overline{(e^{i\theta}\fc)}$ or $(e^{i\theta}\fc)^-=0$ with $\theta$ and $\theta+\pi$ giving the same equations. The points $\theta=\frac{\pi}{4}$ (or $\frac{3\pi}{4}$) are special points. A little algebra shows that at $\theta=\frac{\pi}{4}$, equations (\ref{famKW}) become
\begin{equation}\begin{aligned}\label{KWeq1}F_A-\frac{1}{2}[\phi\wedge\phi]&=*D_A\phi,\\
D_A*\phi&=0.\end{aligned}\end{equation}
This single system is often referred to as ``the'' Kapustin-Witten equation.

There is an alternate form of equation (\ref{famKW}) that will be useful when making estimates.  By adding the self and anti-self dual parts of equation (\ref{famKW}) and rearranging, we get:

\begin{equation}\label{famKWalt}
F_A-\frac{1}{2}[\phi\wedge\phi]=-\cot(2\theta)D_A\phi+\csc(2\theta)*D_A\phi.
\end{equation}

\section{Relationship with the Complex Yang-Mills Equations}\label{relcplxym}

The complex Yang-Mills functional is defined in any dimension as the norm squared of the complex curvature. This reduces to the real Yang-Mills functional when the complex part of the connection vanishes. It is useful to introduce some notation.  First, recall that \begin{equation}\fc=(D_A+i\phi)^2=F_A-\frac{1}{2}[\phi\wedge \phi]+iD_A\phi.\end{equation} The complex Yang-Mills functional is then written as
\begin{equation}\begin{aligned}\ym_\cplx(D_A+i\phi)&=\int_M(|F_A-\frac{1}{2}[\phi\wedge\phi]|^2+|D_A\phi|^2)d\mu\\
&=\int_M\left(-tr\left((F_A-\frac{1}{2}[\phi\wedge\phi])\wedge*(F_A-\frac{1}{2}[\phi\wedge\phi])+[D_A\phi\wedge*D_A\phi]\right)\right)\\
&=\int_M\left(-tr\left((F_A-\frac{1}{2}[\phi\wedge\phi]+iD_A\phi)\wedge*(F_A-\frac{1}{2}[\phi\wedge\phi]-iD_A\phi)\right)\right)\\
&=-\int_Mtr(\fc\wedge*\overline{\fc}).
\end{aligned}\end{equation}

The Euler-Lagrange equations for this functional are
\begin{equation}\label{elym}
\begin{aligned}
D^*_A(F_A-\frac{1}{2}[\phi\wedge\phi])+[D_A\phi\wedge\phi]=&0\\
D^*_AD_A\phi+[(F_A-\frac{1}{2}[\phi\wedge\phi])\wedge\phi]=&0.
\end{aligned}
\end{equation}
These equations are not elliptic, even after the real gauge equivalence is accounted for, so it is necessary to add the moment map condition, which takes into account in some geometric fashion, the action of the complex part of the gauge group:
\begin{equation}\label{momap1}D_A*\phi=0.
\end{equation}

The interaction between the complex gauge action and the functional is not straightforward. To be sure of obtaining solutions, the full Yang-Mills function with a complex gauge term can be treated. We will define this as the augmented complex Yang-Mills functional:
$$A\ym_\cplx(D_A+i\phi)=\int_M|F_A-\frac{1}{2}[\phi\wedge\phi]|^2+|D_A\phi|^2+|D_A*\phi|^2d\mu.$$
The Yang-Mills equations with the moment map condition are clearly critical points of this equation.

The relationships between the Yang-Mills functional and the Kapustin-Witten equation are most easily described by considering the configuration space, $Q_\cplx$, as an affine space based on the complexification of $\nobreak C^\infty(T^*M\otimes\gg)$. Complex conjugation is well defined and naturally extends. The 4-dimensional topological term for the Chern class
$$\int tr(F_A\wedge F_A)$$
extends algebraically to the complex function on $Q_\cplx$:
$$\begin{aligned}\int tr(\fc\wedge\fc)=&\int tr\left((F_A-\frac{1}{2}[\phi\wedge\phi])\wedge(F_A-\frac{1}{2}[\phi\wedge\phi])\right.\\
&\left.-(D_A\phi\wedge D_A\phi)+2i(F_A-\frac{1}{2}[\phi\wedge \phi]\wedge D_A\phi)\right)\end{aligned}.$$
If we algebraically extend the definition of the Chern-Simons functional, we obtain the same relationships as in the real case. We clarify this in the next proposition

In the following, we assume $\ds \ddx{t}\big|_{t=0}D_{A_t}+i\phi_t=B+i\psi$. Then
$$\frac{d}{dt}\fc=[D_A+i\phi,B+i\psi],$$
and
$$\begin{aligned}
\frac{d}{dt}\int tr(\fc\wedge\fc)&=2\int tr\left(\fc\wedge(D_A+i\phi)(B+i\psi)\right)\\
&=2\int tr(\fc\wedge(B+i\psi))-2\int tr\left(D_A+i\phi(\fc)\wedge(B+i\psi)\right)\\
&=2\int_{\partial M^4}tr(\fc\wedge(B+i\psi)).
\end{aligned}
$$
Here, the extension of the Bianchi identity, $D_AF_A=0$, to the complex version $D_A+i\phi(\fc)=0$ is an algebraic computation.

$$\begin{aligned}
\frac{d}{dt}(D_{A_t}+i\phi_t)(\fc)=&[B+i\psi,\fc]+\left(D_{A_t}+i\phi_t\right)[D_{A_t}+i\phi_t,B+i\psi]\\
=&[B+i\psi,\fc]+[\fc,B+i\psi]=0
\end{aligned}$$

We wish to define the complex Chern-Simons functional so that
$$\frac{d}{dt}(CS(D_A+i\phi))=\text{tr}(\fc\cdot(B+i\psi))$$
and get the same relationship as in the real case.
\begin{dfn}$\ds CS(d+A+i\phi)=tr(D(A+i\phi)\wedge(A+i\phi))+\frac{1}{6}tr((A+i\phi)^3).$\end{dfn}
This leads to
\begin{prop}$\int_M tr(\fc\wedge\fc)=8\pi^2k+\int_{\partial M}CS(D_A+i\phi).$\end{prop}
Here, the $8\pi^2k$ is a topological term if $\partial M^4=0$, and defines the ambiguity in the Chern-Simons functional when $\partial M^4=X^3\ne\varnothing$. It is important to note that the ambiguity is only in the real part of the Chern-simons functional on $X^3=\partial M^4$; the complex part is well-defined. The definition we give above is for connections on a trivial bundle over $X^3=\partial M^4$. For bundles which are not trivial, the functional is defined relative to a base point.

Recall that $\ds T^+=\frac{1}{2}(T+*\overline{T})$ and $\ds T^-=\frac{1}{2}(T-*\overline{T})$, so it follows that $*T^+=\overline{T}^+, *T^-=-\overline{T}^-$ and $$Re(T^+\wedge T^-)=-Re\left(tr(T^+T^-)\right)d\mu=<T^+,T^->d\mu=0.$$
Then
\begin{equation}
\begin{aligned}
-|T|^2d\mu&=Re\ tr(T\wedge*\overline{T})=Re\ tr\left((T^++T^-)\wedge(T^+-T^-)\right)\\
&=Re\ tr\left((T^++T^-)\wedge((T^++T^-)-2T^-)\right)\\
&=Re\ tr(T\wedge T)-2|T^-|^2d\mu.
\end{aligned}
\end{equation}
We let
\begin{equation}\begin{aligned}
\ym_\cplx(D_A+i\phi)=&\int_{M^4}|\fc|^2d\mu=\int|e^{i\theta}\fc|^2d\mu\\
=&-Re\int_{M^4}tr(e^{i\theta}\fc\wedge e^{i\theta}\fc)+2\int_{M^4}|(e^{i\theta}\fc)^-|d\mu.
\end{aligned}\end{equation}
Here we have applied the algebra to $T=e^{i\theta}\fc$. Since $\int_M tr(\fc\wedge\fc)=8\pi^2k+\int_{\partial M}CS(D_A+i\phi)$, we add to the previous argument.
$$\begin{aligned}
\ym_\cplx(D_A+i\phi)&=-Re\ e^{2i\theta}\int tr(\fc\wedge\fc)+2\int_{M^4}|(e^{i\theta}\fc)^-|^2d\mu\\
&=-Re\ e^{2i\theta}\left(8\pi^2k+\int_{\partial M^4}CS(A+i\phi)\right)+2\int_{M^4}|(e^{i\theta}\fc)^-|^2d\mu.
\end{aligned}
$$
 Since the other contributions are topological $(\partial M^4=\varnothing)$ or boundary data $(\partial M^4\ne \varnothing)$, we are at a minimum of the functional if $(e^{i\theta}\fc)^-=0$. Hence the Euler-Lagrange equations are satisfied. Such a minimum is referred to as a ``topological'' minimum.

There are two corollaries of the proposition.
\begin{cor}If $P$ is a bundle over a compact manifold with boundary, the only solutions to the Kapustin-Witten equations with $\theta\ne(0,\frac{\pi}{2})$ are flat connections.
\begin{proof} Here we note that if $T^-=0$, then $\ds |T|^2d\mu=<T,T^+-T^->d\mu=-tr(T\wedge T)$. So if $\ds(e^{i\theta}\fc)^-=0$,
\begin{equation}
\begin{aligned}
\int|\fc|^2d\mu&=\int_M|e^{i\theta}\fc|^2d\mu=-\int_Mtr\ e^{2i\theta}(\fc\wedge\fc)\\
&=-e^{2i\theta}8\pi^2k.
\end{aligned}
\end{equation}
This can only be true if $k=0$ and $\fc=0$.
\end{proof}\end{cor}
\begin{cor}Suppose $D_A+i\phi$ is a solution of the Kapustin-Witten equations on $X^3\times[-L,L]$.  then
\begin{equation}e^{2i\theta}\left(\int_{X_{-L}}CS(A+i\phi)-\int_{X_{+L}}CS(A+i\phi)\right)=\int_{X^3\times[-L,L]}|\fc|^2d\mu\end{equation}
\begin{proof}
Here, since $D_A+i\phi\big|_{X^3\times t}$ is a continuous family of connections, the discrete topological term will not appear unless a discrete gauge transformation in $C^\infty(M^3,P|M^3)$ is applied. This is an integration by parts formula.\end{proof}\end{cor}

\section{Weitzenb\^och Formulas}
In order to be able to use the solutions of the Kapustin-Witten equations in a meaningful way in geometry, it is necessary to get some estimates on the behavior of solutions.  Although much of the context in which we hope to use the equations would require compactness of the space of solutions, we should remember that the space of solutions of Hitchin's equations in $2$ dimensions, one of the models we are using for equations, is not compact.  However, associated with the complex connection $d+A+i\phi$, in two dimensions we have a quadratic differential $\varphi=tr(\phi)^2$. In two dimensions, we have the following theorem:
\begin{thm}If $d+A+i\phi$ is a solution of Hitchin's equations
$$F_A-\frac{1}{2}(\phi\wedge\phi)=D_A\phi=D_A^*\phi=0$$
on $\Sigma$, then
\begin{enumerate}[a)]
\item the $(0,2)$ part of $tr(\phi)^2=\varphi$ is a holomorphic quadratic differential on $\Sigma$ and
\item the set of solutions of Hitchin's equations with $\int_\Sigma|\phi|^2d\mu\le K$ is compact.
\end{enumerate}
\end{thm}
Hence $tr(\phi)^2=\varphi$ gives us an auxiliary tool for studying the equations. Of course, $\varphi=tr(\phi)^2$ does form a quadratic differential in all dimensions.  Unfortunately, we do not see how to interpret the equations in this context.

However, we can get some meaningful estimates in all dimensions.  These come from the full complex Yang-Mills equations for the one form $\phi$
\begin{equation*}
D_A^*D_A\phi+\left[F_A-\frac{1}{2}[\phi,\phi],\phi\right]=0,
\end{equation*}
and the moment map equation
\begin{equation*}
D_A^*\phi=0.
\end{equation*}
There is a Weitzenb\^och formula relating the exterior differential with a full differential. It is useful to recall the well-known Bochner formula for one forms on $M$.  We denote the full covarant derivative by $\nabla$ and exterior differentiation by $d$.  For $\alpha$ a one-form on $M$
\begin{equation} \nabla^*\nabla\alpha-\text{Ricci}\ \alpha=(d^*d+dd^*)\alpha.\end{equation}
If we couple this to a connection $A$, we get
\begin{thm} (Weitzenb\^och)
\begin{equation}\nabla_A^*\nabla_A\phi-\text{\emph{Ricci}}\ \phi-[F_A,\phi]=D_A^*D_A\phi+D_AD_A^*\phi.\end{equation}
\end{thm}
\begin{thm}\label{rict}If $D_A+i\phi$ is a solution of the complex Yang-Mills equations, then
\begin{equation}\label{ricf}\nabla_A^*\nabla_A\phi-\text{\emph{Ricci}}\ \phi-\frac{1}{2}[[\phi,\phi],\phi]=0.
\end{equation}
\end{thm}
\begin{cor} \label{ric2}If $\phi$ is the one form for a solution of the complex Yang-Mills equation, then
\begin{equation*}\frac{1}{2}\Delta|\phi|^2-|\nabla_A\phi|^2-|[\phi,\phi]|^2-\text{\emph{Ricci}}(\phi,\phi)\ge0.\end{equation*}
\end{cor}
This formula is proved by taking the inner product of (\ref{ricf}) with $\phi$.

This equation, which applies to solutions of the complex Yang-Mills equations, has three interesting consequences.  By integrating (\ref{ric2}) over $M$ we find:

\begin{cor} If $M$ is a compact manifold with positive Ricci curvature, solutions of the complex Yang-Mills equations reduce to solutions of the real Yang-Mills equations with $\phi\equiv0$.\end{cor}
\begin{cor} \label{bprime}Assume $B'\subset B\subset M$, where the balls $B'$ and $B$ have metric properties close to $B^4\subset B^4_2\subset\real^4$.  Then for solutions of the complex Yang-Mills equations
$$\max_{x\in B'}|\phi^2(x)|\le C\iint_B|\phi|^2d\mu.$$
\end{cor}

The proof is a straightforward application of the maximum principle using (from Corollary \ref{ric2})
$$\frac{1}{2}\Delta(|\phi^2|)\ge-\text{Ricci}(\phi,\phi).$$
The restriction on $B'\subset B$ is applied by taking small balls $B_\rho\subset B\subset M$ and dilating the metric by $\rho\inv$, and making estimates in the dilated metric.  Hence this estimate can be applied on any manifold.

We now use the Kapustin-Witten equations.
\begin{cor}If $D_A+i\phi$ is a solution of the Kapustin-Witten equations for $\theta\in(0,\frac{\pi}{2})$ and $B'\subset B$ as in Corollary (\ref{bprime}), then
\begin{center}$$\begin{aligned}
(a)&\iint_{B'}\left(|\nabla_A\phi|^2+\frac{1}{2}|[\phi,\phi]|^2\right)d\mu\le\overline{C}\iint_B|\phi|^2d\mu\\
(b)&\ds\iint_{B'}|F_A|^2d\mu\le C(\theta)\iint_B|\phi|^2d\mu
\end{aligned}$$\end{center}
\begin{proof}{To obtain (a), multiply (\ref{ric2}) by $f(x)\ge0$, $f(x)\equiv1$ on $B'$ and $f(x)\equiv0$ on $M\setminus B$.  Then
$$\begin{aligned}
\iint_{B'}|\nabla_A\phi|^2+\frac{1}{2}|[\phi,\phi]|^2d\mu&\le\iint_Bf(x)\left(|\nabla_A\phi|^2+\frac{1}{2}|[\phi,\phi]|^2\right)d\mu\\
&\le-\iint_B\text{Ricci}(\phi,\phi)f(x)d\mu+\frac{1}{2}\iint_B\Delta f|\phi|^2d\mu\\
&\le\overline{C}\iint_B|\phi|^2d\mu.
\end{aligned}$$

To obtain (b), we note from (\ref{famKWalt}) that
$$F_A=(F_A-[\phi,\phi])+[\phi,\phi]=-\cot{2\theta}D_A\phi+\csc{2\theta}(*D_A\phi)+[\phi,\phi]$$
immediately implies
$$|F_A|\le|[\phi,\phi]|+2(\cot{2\theta}+\csc{2\theta})|\nabla_A\phi|.$$
The result follows.}
\end{proof}
\end{cor}

To obtain further estimates, it is necessary to make a gauge choice.  We recall the basic gauge-fixing theorem.
\begin{thm}(Uhlenbeck)\cite{springerlink:10.1007/BF01947069} There exists $K>0$, such that if $\nobreak\ds \iint_{B'}|F_A|^2d\mu<K$, then there exists a choice of gauge in which
$$\begin{aligned}
(a)&d^*A=0\\
(b)&\left(\iint_{B'}|A^4|d\mu\right)^{\frac{1}{2}}\le\iint_{B'}|F_A|^2d\mu.
\end{aligned}$$
\end{thm}

This theorem is originally proved in $B_1^4=\{x:|x|\le1\}$.  However, since the metrics are close, it is not difficult to make a gauge change to $\nabla^*A=0$ from $\ds\sum^4_{j=1}\frac{\partial}{\partial x^j}A_j=0$.

Once the gauge is fixed, the Kapustin-Witten equations are manifestly elliptic with a quadratic non-linearity.  This is most easily seen by setting
$$\begin{aligned}
Q_1&=\cos{\theta}A-\sin{\theta}\phi\\
Q_2&=\sin{\theta}A+\cos{\theta}\phi.
\end{aligned}$$
Now
$$\begin{aligned}
d^+Q_1&+\cos{\theta}\frac{1}{2}[A,A]^+-\sin{\theta}[A,\phi]^+=0\\
d^*Q_1&-\sin{\theta}[A,*\phi]=0\\
d^-Q_2&+\sin{\theta}\frac{1}{2}[A,A]^-+\cos{\theta}[A,\phi]^-=0\\
d^*Q_2&+\cos{\theta}[A,*\phi]=0.
\end{aligned}$$
We denoted contraction in the two-form indices of $[A,\phi]$ by $[A,*\phi]$. The quadratic terms in these equations are easily rewritten to be quadratic on $(Q_1,Q_2)$. We now have an equation
$$L(Q)=P(Q,Q)$$
where $L$ is an elliptic operator and $P$ is quadratic.  We now state the main theorem of this section.
\begin{thm}If $\iint_B|\phi|^2<K(\theta)$, where $K(\theta)=K/C(\theta)$, then there exists a gauge in which $(A,\phi)$ and all their derivatives are bounded in interior balls $B''\subset B$.
\end{thm}
This theorem is true for $\theta=0, \pi/2$.  It is simply harder to obtain the estimate for $F_A$.

Since $\phi$ is a one-form, dilation $x\to x/p$ has the effect on $\rho^*\phi=\rho\inv\phi_j(x/p)$
$$\iint_{|x|\le\rho}|\phi|^2=\rho^2\iint|\rho^*\phi|^2d\tilde{\mu}.$$
Any large-scale bound on $\iint_M|\phi|^2d\mu<K$ can be converted into small bounds $\rho^2K=K(\theta)$ by choosing sufficiently small balls.  Hence bounds on $\iint_{|x|\le p}|F_A|^2d\mu$ on a set of solutions imply compactness.

Unfortunately, we see no way to obtain estimates on $\iint_M|\phi|^2d\mu$.  Understanding what happens as $\iint_M|\phi|^2d\mu\to\infty$ is the big open problem in the theory.
\section{The K\"{a}hler Structure on the Space of Complex Connections and the Moment Map for the Action of the Real Gauge Group}\label{sec5}
The integral $\ym$ on the space of complex connections is invariant the real gauge group $\ds P\times_{Ad\ G}G$.  However, the Euler-Lagrange equations are not elliptic.  To remedy this, we
\begin{enumerate}
\item identify a K\"ahler structure on the complex connections,
\item calculate the moment map for the real gauge group, and
\item restrict the equations to the zeros of the moment map.
\end{enumerate}

  The space of complex connections, $\mathbb{Q}_\cplx$, is an affine complex space, and therefore has a flat K\"ahler structure given by
\begin{equation}\begin{aligned}\delta_j(D_A+i\phi)&=B_j+i\psi_j\\
\omega(B_1+i\psi_1,B_2+i\psi_2)&=\int_{M^n}tr(B_1\psi_2-B_2\psi_1)d\mu.\end{aligned}\end{equation}
We can then check that
$$\mu(D_A+i\phi)=D_A^*\phi$$
is the moment map for the gauge action.  If $U$ is in the Lie algebra to the real gauge group, then the infinitesimal change at $D_A+i\phi$ is given by
\begin{equation}[D_A+i\phi,U]=D_AU+i[\phi,U].\end{equation}
The Hamiltonian for the action is
\begin{equation}H(D_A+i\phi,V)=\int_{M^n}tr(D^*_A\phi\cdot V)d\mu\end{equation}
where we identify the dual of the gauge group with itself using the geometry.
$$\begin{aligned}
\delta H(B+i\psi)&=\int_{M^n}tr([B,\phi]\cdot V)d\mu+\int_{M^n}tr(D_A^*\psi\cdot V)d\mu\\
&=\int_{M^n}tr(B\cdot[\phi,V])d\mu-\int_{M^n}tr(\psi\cdot D_AV)d\mu.
\end{aligned}$$
Using the K\"ahler form, we get
$$\begin{aligned}
&\int_{M^n}tr(B\cdot\delta\phi)d\mu-\int_{M^n}tr(\psi,\delta A)d\mu\\
&=\omega(B+i\psi,\delta A+i\delta\phi)
\end{aligned}$$
Hence $\delta A=D_AV$ and $\delta\phi=[\phi,V]$ as claimed.  There is quite a bit of geometry in the equation $D^*_A\phi=0$ as can be seen by the Weizenbach formulas.
\section{Gradient flow for Chern-Simons}

In this section, we use the fact that the complex Chern-Simons functional is complex analytic on $Q=C^\infty(\Lambda T^*M\otimes\gg_\cplx)$. On a complex manifold, we can associate to each analytic functional $f$ a family of real functions $f_\theta=Re( e^{i2\theta}f)$ all with the same critical points. However, the gradient flows are different for different $\theta$.

We actually need only know that
\begin{equation}d_{B+i\psi}{\mathcal CS}(A+i\phi)=\int_{X^3}tr(\fc\wedge(B+i\psi)).\end{equation}

\begin{prop}The gradient flow for $Re(e^{2i\theta}CS)$ on $X^3$ is the Kapustin-Witten equation for $\theta$ on $X^3\times\real$.
\begin{proof}
$$\begin{aligned}
&d_{B+i\psi}Re(e^{2i\theta}{\mathcal CS}(A+i\phi))\\
=&Re\int tr(e^{2i\theta}\fc\wedge (B+i\psi))\\
=&\int tr(*(Re\ e^{2i\theta}\fc)\cdot B)d\mu-tr(*Im(\ e^{2i\theta}\fc)\cdot\psi)d\mu.
\end{aligned}$$

Hence the geometric gradient is
\begin{equation}\label{gf}\begin{aligned}
\frac{\partial A}{\partial t}&=*Re(e^{2i\theta}\fc)\\
\frac{\partial \phi}{\partial t}&=-*Im(e^{2i\theta}\fc).
\end{aligned}\end{equation}

Note that we are extending the flow for real Chern-Simons, which is computed using the geometric gradient. There is no symplectic structure in that theory. However, the symplectic gradient for the complex functional $Re(e^{2i\theta}f)$ is the geometric gradient for $Im(e^{2i\theta}f)=Re(e^{2i(\theta+\frac{3\pi}{4})}f)$ and only involves a change in $\theta$.

In any case, we regard $\frac{\partial A}{\partial t}=F_{t\cdot}$ and $\frac{\partial \phi_\cdot}{\partial t}=(D\phi)_{t\cdot}$ when $A=\sum_j A_jdx^j+0dt$ and $\phi=\sum_j\phi_jdx^j+0dt$. The condition on $A$ is a gauge choice, but the condition on $\phi$ is an assumption.

We now have, using $\cdot$ to indicate a spacial index

\begin{equation}(F_{t\cdot}-[\phi_t\wedge\phi_\cdot])-i(D\phi)_{t\cdot}=e^{2i\theta}*((\fc)_{\cdot\cdot})\end{equation}
or
$$e^{-i\theta}((F_{t\cdot}-[\phi_t\wedge\phi_\cdot])-i(D\phi)_t)=e^{i\theta}*((\fc)_{\cdot\cdot}).$$
This translates into
$$e^{i\theta}\fc=*\overline{e^{i\theta}\fc}$$
which is a succinct way of writing the Kapustin-Witten equations.
\end{proof}
\end{prop}

Up until now, we have been ignoring the moment map condition from section \ref{relcplxym}, $D_A*\phi=0$, which is a purely spacial condition on $Q|(X^3,t)$ rather than a Partial Differential Equation on $X^3\times\real$ because $\phi=\sum_{j=1}^3\phi_jdx^j$. However, we note:
\begin{prop}
The condition $D_A*\phi=0$ is fixed under the gradient flow for $Re(e^{2i\theta}CS(D_A+i\phi))$.
\begin{proof}
To show that $D_A*\phi=0$ we will actually show that $D^*_A\phi=0$ where $D^*_A=*D_A*$. Since $\phi_4=0$ and $\sum_{j=1}^3\phi_jdx^j$, the calculations are on a 3-manifold.
$$\begin{aligned}
\ddx{t}(*D_A*\phi)=&*[\frac{\partial A}{\partial t}\wedge*\phi]+D^*_A\frac{\partial \phi}{\partial t}\\
=&*Re(e^{2i\theta}[*\fc\wedge*\phi])+*D_A*(-Im(e^{2i\theta}*\fc))\\
=&*Re(e^{2i\theta}[\fc\wedge\phi])-*Im(e^{2i\theta}D_A\fc).
\end{aligned}$$
Using the Bianchi identity, we get $D_A\fc=[\fc\wedge i\phi]$. Substituting this in the previous equations gives:
$$\begin{aligned}
&*Re(e^{2i\theta}[\fc\wedge\phi])-*Im(e^{2i\theta}[\fc\wedge i\phi])\\
=&*Re(e^{2i\theta}[\fc\wedge\phi])-*Re(e^{2i\theta}[\fc\wedge \phi])\\
=&0.
\end{aligned}$$
\end{proof}
\end{prop}

\bibliographystyle{plain}
\bibliography{kapwitref}

\end{document}